# DISCUSSION: THE DANTZIG SELECTOR: STATISTICAL ESTIMATION WHEN $p$ IS MUCH LARGER THAN $n$

By Michael P. Friedlander and Michael A. Saunders

*University of British Columbia and Stanford University*

**1. Computational considerations.** When Lasso [11] was proposed, it was a computational challenge to solve the associated quadratic program

$$\text{Lasso}(t) \qquad \min_\beta \tfrac{1}{2}\|y - X\beta\|_2^2 \quad \text{s.t.} \quad \|\beta\|_1 \leq t$$

given just a single parameter $t$. Two active-set methods were described in [11], with some concern about efficiency if $p$ were large, where $X$ is $n \times p$. Later when basis pursuit de-noising (BPDN) was introduced [2], the intention was to deal with $p$ very large and to allow $X$ to be a sparse matrix or a fast operator. A primal–dual interior method was used to solve the associated quadratic program, but it remained a challenge to deal with a single parameter.

The authors' new Dantzig Selector (DS) also assumes a specific parameter. It is helpful to state the BPDN and DS models together:

$$\text{BPDN}(\lambda) \qquad \min_{\beta,r} \lambda\|\beta\|_1 + \tfrac{1}{2}\|r\|_2^2 \quad \text{s.t.} \quad r = y - X\beta,$$

$$\text{DS}(\lambda) \qquad \min_{\beta,r} \|\beta\|_1 \quad \text{s.t.} \quad \|X^T r\|_\infty \leq \lambda, \qquad r = y - X\beta.$$

For reference purposes we also state the corresponding dual problems:

$$\text{BPdual}(\lambda) \qquad \min_r -y^T r + \tfrac{1}{2}\|r\|_2^2 \quad \text{s.t.} \quad \|X^T r\|_\infty \leq \lambda,$$

$$\text{DSdual}(\lambda) \qquad \min_{r,z} -y^T r + \lambda\|z\|_1 \quad \text{s.t.} \quad \|X^T r\|_\infty \leq \lambda, \qquad r = Xz.$$

We congratulate the authors on justifying their Dantzig Selector on detailed statistical grounds while also investigating a primal–dual interior method suitable for a sparse or fast-operator $X$ and making codes available through $\ell_1$-magic [1]. The attraction of a pure linear programming (LP) formulation is understandable. Our aim here is to help explore the prospects for both interior and simplex implementations of DS, and to compare with BPDN.

The vectors $r = y - X\beta$ and $s = -X^T r$ are used often below.

We now know that the Homotopy [5, 7, 8] and LARS [6] algorithms can solve BPDN($\lambda$) for all $\lambda \geq 0$, and their active-set continuation approaches are remarkably efficient if the computed $\beta$ remains sufficiently sparse. Nevertheless, most of our discussion involves a single $\lambda$, and although Lasso came before basis pursuit, we refer mostly to the de-noising problem BPDN($\lambda$) because its $\lambda$ is directly comparable to the DS parameter.

Note from BPdual($\lambda$) that an *optimal* basis-pursuit solution provides a *feasible* solution to DS($\lambda$). Both approaches constrain $\|X^T r\|_\infty \leq \lambda$ while keeping $\|\beta\|_1$ "small," but BPDN strikes a further balance by giving a slightly larger $\|\beta\|_1$ and a slightly smaller $\|r\|_2^2$.

**2. The DS implementation.** The authors eliminate $r$ from DS($\lambda$) and formulate their model as the LP problem

$$\text{(DS)} \quad \min_{\beta, u} \mathbf{1}^T u \quad \text{s.t.} \quad -u \leq \beta \leq u, \quad -\lambda \mathbf{1} \leq X^T(y - X\beta) \leq \lambda \mathbf{1},$$

for which $\ell_1$-magic's MATLAB primal–dual interior solver `l1dantzig_pd` [9] is designed. The main work per iteration lies in solving a $p \times p$ symmetric system

$$(2.1) \quad H\Delta\beta = r_3, \quad H \equiv D_{12} + X^T(XD_{34}X^T)X,$$

where $D_{12}$ and $D_{34}$ are positive definite diagonal matrices. This system is solved in `l1dantzig_pd` using a dense or sparse factorization of $H$ if $X$ is explicit, or the conjugate-gradient method if $X$ is an operator.

To save work when $n \ll p$, the authors suggest reducing (2.1) to an $n \times n$ system that involves the matrix $I + (XD_{34}X^T)(XD_{12}^{-1}X^T)$. Unfortunately this loses symmetry (unnecessarily) and becomes increasingly hazardous as iterations proceed because $D_{12}$ approaches singularity. It is hard to recommend this approach except perhaps for the early iterations.

**3. Test data.** Following $\ell_1$-magic's example, in MATLAB we generated data $X, y$ depending on dimensions $n, p, T$ as follows:

```
rand('state',0);                  % initialize generators
randn('state',0);
q       = randperm(p);            % random +/-1 signal
q       = q(1:T);
beta    = zeros(p,1);
beta(q) = sign(randn(T,1));
[X,R]   = qr(randn(p,n),0);
```

```
X       = X';                     % n x p measurement matrix
y       = X*beta + 0.005*randn(n,1); % noisy observations
```

Thus, $X$ is dense with orthogonal rows ($XX^T = I$) and $\beta$ should have $T$ components close to $\pm 1$. We used $\lambda = $ `3e-3` for all test cases. Times are cpu seconds on a 3.2 GHz Linux Intel Pentium 4 with 2 GB of memory.

**4. DS and BPDN with interior solvers.** To compare with more general primal–dual interior solvers, we considered two formulations of the DS problem and also the BPDN formulation in [3]:

(DS1)
$$\min_{v,w,s} \mathbf{1}^T(v+w)$$
$$\text{s.t.} \quad [X^TX \ -X^TX \ I] \begin{bmatrix} v \\ w \\ s \end{bmatrix} = X^Ty, \qquad v, w \geq 0, \|s\|_\infty \leq \lambda,$$

(DS2)
$$\min_{v,w,r,s} \mathbf{1}^T(v+w)$$
$$\text{s.t.} \quad \begin{bmatrix} X & -X & I & \\ & & X^T & I \end{bmatrix} \begin{bmatrix} v \\ w \\ r \\ s \end{bmatrix} = \begin{bmatrix} y \\ 0 \end{bmatrix}, \qquad v, w \geq 0, \|s\|_\infty \leq \lambda,$$

(DS3)
$$\min_{v,w,r} \lambda \mathbf{1}^T(v+w) + \tfrac{1}{2} r^T r$$
$$\text{s.t.} \quad [X \ -X \ I] \begin{bmatrix} v \\ w \\ r \end{bmatrix} = y, \qquad v, w \geq 0,$$

TABLE 1
*Dense orthogonal X*

| Sizes | | | (DS) | (DS1) | (DS2) | | (BPDN) | | |
|---|---|---|---|---|---|---|---|---|---|
| $n$ | $p$ | $T$ | l1magic | Pdco | Pdco | Cplex | Pdco | Cplex | Greedy |
| 120 | 512 | 20 | 1.2 | 2.7 | 3.9 | 1.7 | 0.2 | 0.5 | 0.1 |
| 240 | 1024 | 40 | 6.9 | 16.1 | 24.5 | 16.5 | 1.0 | 4.9 | 0.2 |
| 360 | 1536 | 60 | 20.8 | 48.9 | 75.6 | 58.1 | 2.4 | 15.3 | 0.4 |
| 480 | 2048 | 80 | 46.9 | 110.7 | 171.3 | 122.8 | 5.0 | 34.3 | 1.0 |
| 720 | 3072 | 120 | 149.4 | 349.7 | 550.1 | 391.6 | 14.6 | 109.6 | 3.4 |
| 960 | 4096 | 160 | 349.1 | 814.0 | 1275.7 | 855.4 | 31.8 | 245.3 | 9.2 |

Cpu time for 15 iterations of three primal–dual interior solvers and $T$ iterations of a Homotopy/LARS-type greedy algorithm.

where $\beta = v - w$ and $\|\beta\|_1 = \mathbf{1}^T(v+w)$, and we expect few nonzero elements in $v$ and $w$. We applied l1dantzig_pd [9], Pdco [10] and the Cplex barrier LP/QP solver [4] to the relevant problem formulations (see Table 1). With $X$ dense, all solvers use dense Cholesky factors of matrices of the form $H = AD_1A^T + D_2$, where $A$ denotes the corresponding constraint matrix and $D_1, D_2$ are positive diagonal matrices that change each iteration. (We modified l1dantzig_pd slightly to ensure that its $H$ was recognized to be symmetric positive definite.)

Table 1 shows computation times on increasingly large problems. The $\ell_1$-magic solver is specialized to problem (DS) and operates with $X^TX$ only once, whereas Pdco must double-handle that matrix in (DS1) and has $n+p$ general constraints to deal with in (DS2). Cplex barrier solves all (DS2) examples in times midway between those for the other two solvers.

We see that the solution times are rather large for all DS formulations and solvers. In contrast, Pdco is quite efficient on the BPDN problems, primarily because there are only $n$ general constraints. A minor specialization to avoid double-handling $X$ would reduce times further. We expected the Cplex barrier QP solver to perform comparably on the BPDN examples (since its interior algorithm is similar to that in Pdco). In case unbounded variables were not handled well by Cplex's barrier implementation, we added bounds on $r$ enforcing $\|r\|_\infty \leq \|y\|_2$, but the times remained essentially the same.

The greedy method listed in Table 1 is an experimental MATLAB active-set method intended for problem BPDN($\lambda$) with a specific $\lambda$. Like Homotopy and LARS, it starts with $\beta = 0$ and selects one parameter at a time—in this case, the one whose dual constraint is most violated. It required exactly $T$ iterations on these examples, each involving multiplications with $X$ and $X^T$ (to compute $r$ and $s$) and a QR factorization of $S$, the columns of $X$ chosen so far.

If $T$ were changed in each test case, the solution times for the interior methods would be essentially unaltered, but for the greedy method they would change in proportion to $T$.

If $X$ is sparse but $X^TX$ is not, interior solvers on (DS2) could potentially be more efficient than on (DS) or (DS1). However, in trying to generate random sparse examples we found that the expected $T$ nonzero parameters were not correctly identified. The sparse $X$ case remains for study. Both l1dantzig_pd and Pdco allow $X$ to be an operator, but we have not compared those options.

Donoho and Tsaig [5] give related computational results for Homotopy, Pdco and simplex for the basis-pursuit case $r = 0$ (another LP setting!), with $n, p, T$ as large as $1600, 4000, 320$ and dense $X$ drawn from the Uniform Spherical Ensemble. Again the greedy Homotopy approach performs best.

**5. DS and the simplex method.** It seems clear that formulations (DS) and (DS1) are not well suited to general-purpose simplex codes for two reasons: the presence of a potentially dense $X^T X$, and the large number of constraints (viz., $p$).

For a time, we thought that formulation (DS2) might be *ideal* for large-scale simplex solvers such as in Cplex. This would be for a specific $\lambda$ and values of $T$ up to a few hundred, or a few thousand if $X$ were sparse. If the initial basis includes $r$ and $s$ (with nonbasic variables $v = w = 0$), the initial primal and dual variables can be cheaply computed from

$$(5.1) \quad \begin{bmatrix} I & \\ X^T & I \end{bmatrix} \begin{bmatrix} r \\ s \end{bmatrix} = \begin{bmatrix} y \\ 0 \end{bmatrix} \quad \text{and} \quad \begin{bmatrix} I & X \\ & I \end{bmatrix} \begin{bmatrix} \bar{r} \\ \bar{s} \end{bmatrix} = \begin{bmatrix} 0 \\ 0 \end{bmatrix}.$$

Note that the initial dual values $\bar{r} = \bar{s} = 0$ are dual feasible, and $r$ will remain in the basis throughout. We hoped that the dual simplex method would proceed in an essentially greedy fashion until $T$ components of $v$ or $w$ replaced $T$ components of $s$. The basis would remain almost triangular and therefore easy for a typical sparse LU factorization. If we partition $X = [Z \ S]$ to match the current zero and nonzero parameters, the basis LU factors take the form

$$B \equiv \begin{bmatrix} I & & \bar{S} \\ Z^T & I & \\ S^T & & I \end{bmatrix}$$

$$= \begin{bmatrix} I & & \\ Z^T & I & \\ S^T & & \bar{L} \end{bmatrix} \begin{bmatrix} I & & \bar{S} \\ & I & -Z^T \bar{S} \\ & & \bar{U} \end{bmatrix} \quad \text{with} \quad \bar{L}\bar{U} = -S^T \bar{S},$$

where $\bar{S}$ is the same as $S$ with columns scaled by $\pm 1$ according to whether an element of $v$ or $w$ is basic. The work per iteration with such factors is much the same as for Homotopy/LARS: multiplications by $X$ and $X^T$ and factorization of $S^T S$. (A specialized basis factorization could account for the special structure of $S^T \bar{S}$ and compute a QR factorization of $S$.)

A specialized simplex solver could be constructed to use the same LU factorization even if $X$ is an operator. Ideally, $S$ would be kept in memory as its columns come and go.

Further, we note that if all values of $\lambda$ are of interest, problem (DS2) may be treated as an LP problem with parametric bounds. A simplex-type algorithm for such problems is known [12] that works directly with the original variables and constraints. Thus a Homotopy/LARS-type algorithm does indeed seem practical at first sight.

The most effective Cplex simplex options we could find were *dual simplex*, *no scaling*, *no presolve* and *steepest-edge pricing*. Results are summarized in Table 2. Unfortunately, it appears that simplex methods work in a "far from greedy" manner. In a genuinely optimal solution, many more than $T$

Table 2
*Dense orthogonal X*

| Sizes | | | tol = 0.1 | | | tol = 0.01 | | | tol = 0.001 | | |
|---|---|---|---|---|---|---|---|---|---|---|---|
| $n$ | $p$ | $T$ | itns | $\|S\|$ | Time | itns | $\|S\|$ | Time | itns | $\|S\|$ | Time |
| 120 | 512 | 20 | 20 | 20 | 0.1 | 20 | 20 | 0.1 | 86 | 63 | 0.2 |
| 240 | 1024 | 40 | 58 | 56 | 0.4 | 67 | 57 | 0.4 | 405 | 150 | 2.3 |
| 360 | 1536 | 60 | 187 | 134 | 2.3 | 655 | 156 | 7.8 | 1231 | 215 | 15.1 |
| 480 | 2048 | 80 | 163 | 122 | 3.4 | 549 | 211 | 11.1 | 1277 | 275 | 26.7 |
| 720 | 3072 | 120 | 356 | 223 | 15.3 | 1414 | 317 | 65.0 | 3006 | 420 | 146.6 |
| 960 | 4096 | 160 | 965 | 414 | 80.2 | 6226 | 488 | 574.9 | 9229 | 567 | 891.6 |

Cplex dual simplex on problem (DS2) with loose and tighter termination tolerances.

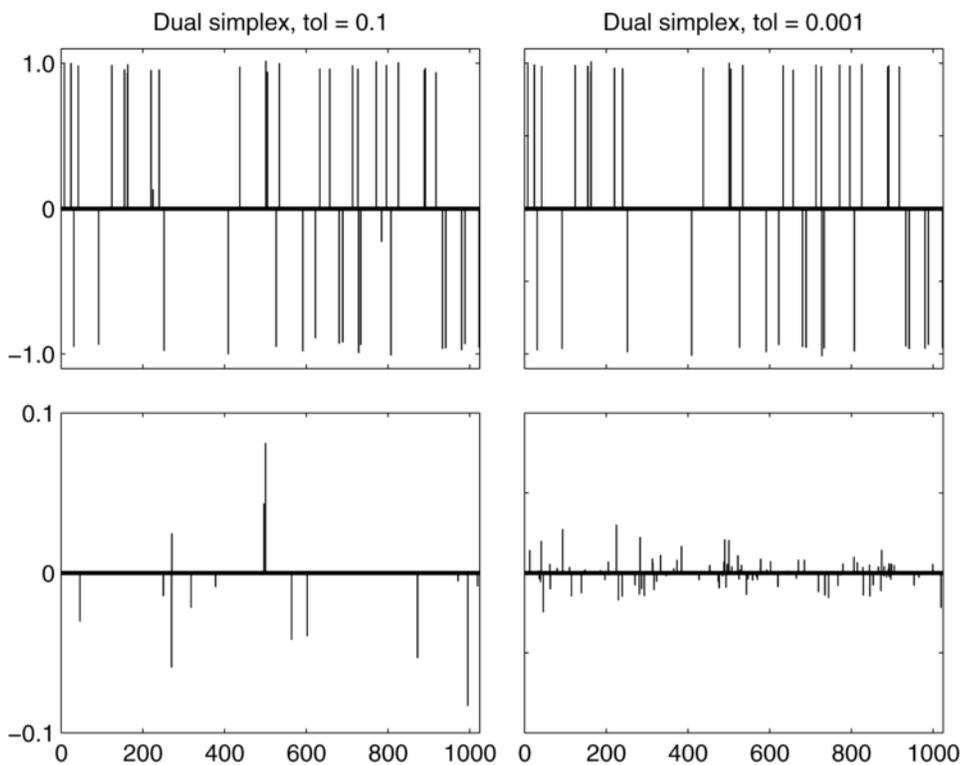

FIG. 1. *Cplex dual simplex method on $240 \times 1024$ problem (DS2) with $T = 40$ nonzero "true" parameter values of $\pm 1$. Plot of significant (top) and small (bottom) solution values with two termination tolerances. More small values imply more simplex iterations and more time per iteration.*

parameters enter the basis, and the number of simplex iterations exceeds $T$ by a huge factor. (Thus, many parameters must be getting selected and then rejected.) This has dampened our optimism for the effectiveness of simplex on large-scale DS problems.

On the other hand, the degree of optimality required can have a profound effect. Table 2 shows the trend with several feasibility and optimality tolerances (tol = 0.1, 0.01 and 0.001). We would normally regard tol = 0.1 as unusually "loose," but Figure 1 emphasizes the benefit of terminating early (at the risk of violating $\|X^T r\|_\infty \leq \lambda$ by as much as tol!).

**6. Conclusions.** We have tested interior solvers on three DS formulations, and compared with three BPDN solvers on the same data. Table 1 results confirm that the larger DS constraint matrix is likely to invoke a high computational cost compared to the Lasso/BPDN model.

In keeping with the DS name, we have also tested some simplex codes (which seem necessary if a range of $\lambda$ values is of interest). Table 2 again predicts a high cost, except perhaps if low-accuracy solutions are acceptable.

Tables 1 and 2 and Figure 1 can be reproduced using the MATLAB scripts in www.cs.ubc.ca/labs/scl/ds_discussion.html.

We emphasize that the solvers tested are *general purpose.* They would all be "happier" if the dense data $X$ were sparse, and none of them takes advantage of the property $XX^T = I$ (which may arise in certain situations). We wish the authors much success in exploring the virtues of their linear DS model for an increasing range of real-world applications.

## REFERENCES

[1] CANDÈS, E. (2007). $\ell_1$-magic. Available at www.l1-magic.org.
[2] CHEN, S. S., DONOHO, D. L. and SAUNDERS, M. A. (1998). Atomic decomposition by basis pursuit. *SIAM J. Sci. Comput.* **20** 33–61. MR1639094
[3] CHEN, S. S., DONOHO, D. L. and SAUNDERS, M. A. (2001). Atomic decomposition by basis pursuit. *SIAM Rev.* **43** 129–159. MR1854649
[4] ILOG CPLEX (2007). Mathematical programming system. Available at www.cplex.com.
[5] DONOHO, D. L. and TSAIG, Y. (2006). Fast solution of $\ell_1$-norm minimization problems when the solution may be sparse. Available at www.dsp.ece.rice.edu/cs/FastL1.pdf.
[6] EFRON, B., HASTIE, T., JOHNSTONE, I. and TIBSHIRANI, R. (2004). Least angle regression (with discussion). *Ann. Statist.* **32** 407–499. MR2060166
[7] OSBORNE, M. R., PRESNELL, B. and TURLACH, B. A. (2000). On the LASSO and its dual. *J. Comput. Graph. Statist.* **9** 319–337. MR1822089
[8] OSBORNE, M. R., PRESNELL, B. and TURLACH, B. A. (2000). A new approach to variable selection in least squares problems. *IMA J. Numer. Anal.* **20** 389–403. MR1773265
[9] ROMBERG, J. (2005). l1dantzig_pd.m. MATLAB solver for DS problem. Available at www.l1-magic.org.

[10] SAUNDERS, M. A. (2005). PDCO. MATLAB software for convex optimization. Available at www.stanford.edu/group/SOL/software/pdco.html.

[11] TIBSHIRANI, R. (1996). Regression shrinkage and selection via the lasso. *J. Roy. Statist. Soc. Ser. B* **58** 267–288. MR1379242

[12] TOMLIN, J. A. (1975). A parametric bounding method for finding a minimum $\ell_\infty$-norm solution to a system of equations. Report SOL 75-12, Dept. Operations Research, Stanford Univ.

DEPARTMENT OF COMPUTER SCIENCE  
UNIVERSITY OF BRITISH COLUMBIA  
VANCOUVER, BRITISH COLUMBIA  
CANADA V6K 2C6  
E-MAIL: mpf@cs.ubc.ca

DEPARTMENT OF MANAGEMENT SCIENCE  
AND ENGINEERING  
STANFORD UNIVERSITY  
STANFORD, CALIFORNIA 94305  
USA  
E-MAIL: saunders@stanford.edu